\documentclass[11pt]{article}
\usepackage{amsmath,amstext,amsthm,amsfonts,amscd,amssymb}
\usepackage[dvips]{graphicx}
\usepackage{epic,eepic}
\usepackage[latin1]{inputenc}
\theoremstyle{plain}
\newtheorem{theorem}{Theorem }[section]
\newtheorem{proposition}[theorem]{Proposition}

\newtheorem{maintheorem}{Theorem}

\theoremstyle{definition}
\newtheorem{remark}[theorem]{Remark}

\newtheorem{definition}[theorem]{Definition}

\newtheorem*{ack}{Acknowledgement}

\newcommand{\field}[1]{\mathbb{#1}}
\newcommand{\real}{\field{R}}
\newcommand{\complex}{\field{C}}
\newcommand{\integer}{\field{Z}}

\newcommand{\vep}{\varepsilon}

\newcommand{\SC}{{\cal C}}

\begin{document}

\title{\Large{\textbf{ Decidability of chaos for some families of
dynamical systems}}}

\author{ Alexander Arbieto\footnote{A. Arbieto is supported by CNPq/Brazil}
\text{ and}
Carlos Matheus\footnote{C. Matheus is supported by Cnpq/Brazil and Faperj/Brazil} }
\date{December 4, 2002}

\maketitle

\begin{abstract} We show that existence of positive Lyapounov exponents and/or 
SRB measures are undecidable (in the algorithmic sense) properties within some parametrized families of interesting
dynamical systems: quadratic family and H\'enon maps. Because the existence of positive exponents (or SRB
measures) is, in a natural way, a manifestation of ``chaos'', these results may be understood as saying that the chaotic character of a dynamical system is undecidable. Our investigation is directly motivated by questions asked by Carleson and Smale in this direction.
\end{abstract}

\section{Introduction}

The notion of decidability goes as far back as to the third decade of the previous century. At the time even before computers were actually built, an important issue was which functions can be computed by an automatic machine. Notions such as that of recursive functions and Turing machines were proposed to answer such a question. It turned out that all purposed definition of computable functions were equivalent (Church`s thesis). Then a set is called decidable if its characteristic function is computable. Problems of decidability have played an important role in several branches of Mathematics since then (e.g. the word problem, solved by Novikov~\cite{Nov} or the classification problem for 4-manifolds~\cite{Mar}).

Recently, it has been asked by Carleson in his talk on the ICM-1990 (see p. 1246 of~\cite{Ca}) and by Smale~\cite{Smale} (included in the fourteen problem on Smale's paper) whether the set of chaotic systems is decidable, either in a general context or within special relevant families of systems. Our main result, Theorem~\ref{t.um}, is an answer to this question in the setting of certain families of low dimensional maps: the quadratic family and the H\'enon family, where ``chaos'' in this setting will be positive Lyapounov exponents, SRB measures or positive entropy.

\begin{definition} We denote by $Q_a (x) = 1- ax^2, a,x\in\real$ the
\textit{quadratic family} and $H_{a,b}(x,y)=(1-ax^2+y, bx), (x,y)\in\real^2,
a,b\in\real$ the
\textit{H\'enon maps}. 
\end{definition}

We deal with some sets associated to the maps defined above (for the precise definitions see section~\ref{Lemmas}):

\begin{itemize}
\item $exp^+(Q):=\{ a\in (0,2] \ : \ Q_a \textrm{ has positive Lyapounov exponent}\}$ and \\ $exp^+(H):= \{ (a,b)\in (0,2]\times(-b_0,b_0) \ : \ H_{a,b}  \textrm{ has positive Lyapounov}$ $\textrm{exponents}\}$ (with $b_0>0$ is small). 
\item  $srb(Q):=\{ a\in (0,2] \ : \ Q_a \textrm{ has a SRB measure}\}$ and $srb(H):= \{ (a,b)\in (0,2]\times(-b_0,b_0) \ : \ H_{a,b}  \textrm{ has a SRB measure}\}$ (with $b_0>0$ is small). 

\item $ent^+(Q):=\{ a\in (0,2] \ : \ Q_a \textrm{ has positive topological entropy}\}$ 

\end{itemize}

Then our main result is:

\begin{maintheorem}
\label{t.um}

The sets $exp^+(Q)$, $exp^+(H)$, $srb(Q)$ and $srb(H)$ are undecidable. The set $ent^+(Q)$ is decidable.

\end{maintheorem}

We define the technical notions of positive Lyapounov exponents, SRB measures,
topological entropy and decidable sets in section~\ref{Lemmas}. For instance, we 
make some comments
on the problem of decidability of sets associated to dynamical systems.
In~\cite{BSS}, Blum, Shub and Smale present a model for computation over
arbitrary ordered rings $R$ (for example, $R=\real$). For this general setting,
they obtain a theory which reflects the classical theory over $\integer$.
However, one virtue of this theory is it forces the use 
of more algebraic methods,
i.e., classical mathematics, than the approach from logic. In
particular,~\cite{BSS} shows that most Julia sets (a set associated to the
rational maps $g:\overline{\complex}\rightarrow\overline{\complex}$ of the
Riemann sphere $\overline{\complex}$, just like the Mandelbrot set) are undecidable. Next, there is an 
important problem of
decidability of sets associated to dynamical systems, namely : ``Is the
Mandelbrot set decidable ?''. This problem was first proposed by
Penrose~\cite{Pen}, but since the classical computational theory deals with Turing machines (integers inputs), the question makes no sense. Now, with the modern theory cited above, Blum and Smale formalized and solved the question. For references 
and comments on this problem
see~\cite{BS}. We observe that a result by Shishikura~\cite{Shi} says that the
boundary of the Mandelbrot set has Hausdorff dimension equal to $2$. In
particular, the argument given in~\cite{BS} implies that the Mandelbrot is not
decidable. For sake of completeness, we briefly recall some steps of
Blum-Smale's argument. If the Mandelbrot set $M$ is the halting set of some
machine, the boundary of $M$ is the union of sets with Hausdorff dimension at
most $1$. In particular, $\partial M$ has Hausdorff dimension at most $1$, a
contradiction with Shishikura's theorem. These comments shows that the problem 
of decidability of sets associated to dynamical
systems is relevant. 

The paper is organized as follows. In next section we define the
fundamental concepts and state the results used in section~\ref{Proofs}. After that, we give the proof of the
theorem~\ref{t.um}. Finally,
in the last section, we finish this paper with some final remarks on some open
question related with the subject treated here.

\section{Basic Results}\label{Lemmas}
This section is devoted to introduce some basic concepts used in the proof of
our theorem.

{\bf Dynamical Results:}

All of the concepts present here has general statements which holds for any dynamical system. Here we will give adapted definitions for the systems considered in this paper, which are more simple.

\smallskip

a) {\it Lyapounov Exponents}:

\smallskip

The Lyapounov exponent for $Q_a(x)$ and $H_{a,b}(x)$ are defined by: 
$$\lambda(Q_a)=\limsup\limits_{n\rightarrow\infty}\frac{1}{n}\log |DQ_a^n(1)|$$

$$\lambda(H_{a,b})=\limsup\limits_{n\rightarrow\infty}\frac{1}{n}\log \|DH_{a,b}^n(z_0)(0,1)\|$$Here $z_0$ is a critical point in the sense of~\cite{BC2}.

We also remark that if we take $b=0$ in the H\'enon map then the positivity of the Lyapounov exponents of $H_{a,0}$ implies the positivity of the Lyapounov exponent of $Q_a$.

In general we define the Lyapounov exponents for an invariant measure $\mu$ of a  differentiable transformation $T:M\rightarrow M$ like the possible values ($x \ \mu$-a.e.) of: 
$$\limsup\limits_{n\to \infty} \frac{1}{n}\log \|DT^n(x)v\| \textrm{  , }v\in T_xM-\{0\}.$$

Here we have the remarkable result form Benedicks and Carleson~\cite{BC1} (see also~\cite{dMvS} for other proof and some other references):

\begin{theorem}[Benedicks and Carleson]\label{lemma_BC1} For a subset of positive
measure $E\subset (0,2]$, if $a\in E$, then $Q_a$ has positive Lyapounov
exponent.
\end{theorem}

\smallskip

b) {\it Sinai-Ruelle-Bowen Measures}:

\smallskip

A SRB measure for $Q_a(x)$ is an absolutely continuous (with respect to Lebesgue measure) invariant measure for $Q_a(x)$. A measure $\mu$ is a SRB measure for $H_{a,b}(x)$ if there is a positive Lyapounov exponent $\mu$-a.e and the conditional measures of $\mu$ on unstable manifolds are absolutely continuous with respect to the Lebesgue measure induced on these manifolds (for details see~\cite{BY} or~\cite{BC2}).   

Now we have the following result by Jakobson~\cite{J} (see~\cite{dMvS} for an alternative proof):

\begin{theorem}[Jakobson]\label{lemma_J} For a subset of positive Lebesgue measure 
$\SC\subset
(0,2]$, if $a\in\SC$ then $Q_a$ has a SRB measure with positive entropy.
\end{theorem}

c) {\it Topological Entropy}:

\smallskip
 
Let $X$ a compact metric space, $f:X\rightarrow X$ a continuous, $n\in \mathbb{N}$ and $\vep>0$, a subset $E\subset X$ is $(n,\vep)$-separated if for any $x\neq y$ inside of $E$, there is $0\leq j\leq n$ such that $d(f^j(x),f^j(y))>\vep$. We denote $s_n(\vep)$ the largest cardinality of any set $E\subset X$ which is $(n,\vep)$-separated. Then we define:

$$h_{top}(f)=\lim_{\vep \to 0}\limsup\limits_{n\to \infty}\frac{1}{n}\log s_n(\vep)$$

We refer the reader to~\cite{dMvS} for the proof and references of the next result.

\begin{theorem}\label{monot.} For the quadratic family $Q_a$, the map
$a\rightarrow h_{top}(Q_a)$ is monotone.
\end{theorem}

d) {\it Attracting cycles}:

\smallskip

An attracting cycle is a periodic point $p$ for $f$ (of period $n$) such that $|(f^n)'(p)|<1$. Related to this notion we have the following theorem due to Graczyk and \'Swiatek~\cite{GS} and Lyubich~\cite{L}.
 
\begin{theorem}[Graczyk, \'Swiatek and Lyubich]\label{lemma_GS}For an open and dense set of parameters $a\in
(0,4]$ the mapping $Q_a$ has an attracting cycle.
\end{theorem}

\smallskip

{\bf Computational Results:}

\smallskip

Here we give some notions from computation and complexity theory. We invite the reader to see the work of Blum, Shub and Smale~\cite{BSS} for detailed definitions.

We deal with machines over an ordered ring $R$ (e.g. $R=\real$) introduced by~\cite{BSS}. Following their terminology we will call $S\subset R^n$ a \emph{halting set} if there exist a machine $M$ such that for any input $a\in S$ the machine halts in finite time. We also remark that the halting sets are the R.E. sets of their recursive function theory.

We say that $\Omega\subset R^n$ is a \emph{decidable set} if $\Omega$ and its complement are halting sets. In this case there is a machine which decides for each $z\in R^n$ ``Is $z\in \Omega$?''.

A \emph{basic semi-algebraic} set is a subset of $R^n$ defined by a set of
polynomial inequalities of the form $h_i(x)<0, i=1,\dots,l$, $h_j(x)\leq 0,
j=l+1,\dots, m$.

The following key proposition relates halting sets and basic semi-algebraic sets:

\begin{proposition}\label{prop_BSS} Any R.E. set over $\real$ is the countable 
union of basic
semi-algebraic sets. In particular, any R.E. set over $\real$ has a countable
number of connected components. 
\end{proposition} 

The proof of this proposition is contained in~\cite{BSS}. In fact, they show
that any R.E. set is the countable union of basic semi-algebraic sets. Since a 
basic semi-algebraic set has a finite number of connected components, any R.E.
set has a countable number of connected components. In the same paper, they use this proposition to prove that most Julia sets are undecidable. Here we have used similar arguments, in particular this proposition.

\section{Proof of the theorem}\label{Proofs}

\begin{proof}[Proof of theorem~\ref{t.um}] First, we consider the set
$exp^+(Q)$ of parameters of the quadratic family with positive Lyapounov 
exponents. By a remarkable
result of Benedicks and Carleson (see theorem~\ref{lemma_BC1} below), $exp^+(Q)$ has 
positive Lebesgue measure. However, if $a\in exp^+(Q)$, then $Q_a$ does not
have attracting periodic orbits (see~\cite{dMvS}) and, by theorem~\ref{lemma_GS}
below, the set of the parameters $a$ such that $Q_a$ admits an attracting
periodic orbit (i.e., $Q_a$ is hyperbolic) is an open dense set. In particular,
$int(exp^+(Q))=\emptyset$. Thus, $exp^+(Q)$ has an uncountable number of
connected components. But, proposition~\ref{prop_BSS} below says that any
decidable set has a countable number of connected components. So, $exp^+(Q)$ is
not a decidable set.

Second, we consider the set $exp^+(H)$. Suppose that $exp^+(H)$ is decidable and
denote by $M'$ the machine whose ``halting'' set (see section~\ref{Lemmas}) is $exp^+(H)$. We can define a
new machine $M$ as follows. The inputs are parameters $a\in (0,2]$. Given a such $a$,
we define the input $(a,0)$ for the machine $M'$ and, by hypothesis, $M'$
outputs either $1$ (\emph{yes}) if $(a,0)\in exp^+(H)$ or $0$ (\emph{no}) if
$(a,0)\notin exp^+(H)$. Since, for $b=0$,
$H_{a,b}(x,0)=H_{a,0}(x,0)=(1-ax^2, 0)$, the machine $M$, in fact,
``decides'' $exp^+(Q)$, a contradiction with the previous paragraph. This
concludes the proof for these sets.   

By a Jakobson`s theorem (see theorem~\ref{lemma_J} below), $srb(Q)$
has positive Lebesgue measure. Since the existence of absolutely continuous
invariant measures is an obstruction to the hyperbolicity (i.e., attracting
periodic orbits), $srb(Q)$ has empty interior. In particular, it is sufficient
to make a \emph{mutatis mutandis} argument to conclude the proof for the sets $srb(Q)$ and $srb(H)$.

From the definition of decidable
sets, it is easy that intervals are decidable. Now, theorem~\ref{monot.} says that
the topological entropy of the quadratic family is monotone. In particular,
$ent^+(Q)$ is an interval. This concludes the proof of the theorem. 
\end{proof}

\begin{remark}
In fact, if there exists a periodic orbit for $Q_a(x)$ with period $s$ and $s$ is not of the form $s=2^n$, then the topological entropy is positive. Therefore, following the period doubling cascade we can find the number $c=\sup\{a \ :\ \textrm{$Q_a$ has zero topological entropy}\}.$
\end{remark}

\section{Some Questions}

The fact that positive topological entropy is a decidable property in one dimension, does not permit us to directly reduce the two-dimensional case for the previous case. So, the following very interesting question remains open.

{\bf Question 1:} Is the set $\{h_{top}(H_{a,b})>0\}$ undecidable?

In the same line we have the following question due to Milnor~\cite{Mil}:

{\bf Question 2:} Given an explicit dynamical system (e.g. quadratic family or H\'enon map) and $\vep>0$ is it possible to compute the topological entropy with a maximum error of $\vep$?

An instructive model to explore the complexity of the H\'enon map can be a piecewise linear (possible discontinuous) map on the interval with two parameters. In this case the topological entropy may be calculated using the approach by Widodo~\cite{W}.

\begin{ack}
We are indebted to Marcelo Viana for proposing the problem to us. Also, we are thankful for his encouragement, advices and suggestions. Finally, the authors are grateful to IMPA and his staff. 
\end{ack}

\vspace{1cm}

\noindent    \textbf{Alexander Arbieto} ( alexande{\@@}impa.br )\\
             \textbf{Carlos Matheus} ( matheus{\@@}impa.br )\\
             IMPA, Est. D. Castorina 110, Jardim Bot\^anico, 22460-320 \\
             Rio de Janeiro, RJ, Brazil


\begin{thebibliography}{15}

\bibitem{BC1}
M.~Benedicks and L.~Carleson,
\newblock On iterations of $1-ax^2$ on $(-1,1)$,
\newblock {\em Annals of Math.}, vol.122, 1--25, 1985

\bibitem{BC2}
M.~Benedicks and L.~Carleson,
\newblock The dynamics of the H\'enon map,
\newblock {\em Annals of Math.}, vol.122, 1--25, 1985

\bibitem{BY}
M.~Benedicks and L.S.~Young,
\newblock Sinai-Bowen-Ruelle measures for certain H\'enon maps.
\newblock {\em Invent. Math.}, vol. 112, no 3, 541--576, 1993.


\bibitem{BSS}
L.~Blum, M.~Shub and S.~Smale,
\newblock On a theory of computation and complexity over the real numbers:
NP-completeness, recursive functions and universal machines,
\newblock {\em Bull. of the A.M.S.}, vol.21, 1--46, 1989.

\bibitem{BS}
L.~Blum and S.~Smale,
\newblock The Godel incompleteness theorem and decidability over a ring,
\newblock {\em From Topology to Computation: Proc. of the Smalefest}, 321--339,
1993.

\bibitem{Ca}
L.~Carleson,
\newblock The dynamics of non-uniformly hyperbolic systems in two variables, 
\newblock {\em Proceedings of the International Congress of Mathematicians}, Vol. I, II (Kyoto, 1990),  1241--1247, Math. Soc. Japan, Tokyo, 1991.

\bibitem{dMvS}
W.~de Melo and S.~van Strien,
\newblock One-dimensional dynamics,
\newblock Springer-Verlag, 1993.

\bibitem{GS} 
J.~Graczyk and G.~\'Swiatek,
\newblock Generic hyperbolicity in the logistic family,
\newblock {\em Annals of Math.}, vol.146, 1--52, 1997. 

\bibitem{J}
M.~Jakobson,
\newblock Absolutely continuous invariant measures for one-parameter families of
one-dimensional maps,
\newblock {\em Comm. Math. Phys.}, vol.81, 39--88, 1981.

\bibitem{L}
M.~Lyubich,
\newblock Dynamics of quadratic polynomials. I, II,
\newblock{\em. Acta Math.}, vol.178, no. 2, 185--247, 247--297, 1997. 

\bibitem{Mar}
A.A.~Markov,
\newblock 
\newblock {\em Proceedings of ICM}, 300-306, 1958.

\bibitem{Mil}
J.~Milnor,
\newblock A problem: Is entropy effectively computable?
\newblock {\em http://www.math.sunysb.edu/\~{ }jack/PREPR/index.html}

\bibitem{Nov}
P.S.~Novikov,
\newblock On the algorithmic unsolvability of the word problem in group theory,
\newblock {\em Amer. Math. Soc. Translations}, ser. 2, v. 9, 1--122, 1958.


\bibitem{Pen}
R.~Penrose,
\newblock The Emperor's new mind,
\newblock Oxford Univ. Press, Oxford, 1989. 

\bibitem{Smale}
S.~Smale,
\newblock Mathematical Problems for the next century,
\newblock {\em Mathematics: Frontiers and Perspectives, A.M.S.}, 271--294, 2000.

\bibitem{Shi}
M.~Shishikura,
\newblock The Hausdorff dimension of the boundary of the Mandelbrot set and
Julia set,
\newblock {\em Annals of Math.}, vol.147, 225--267, 1998.

\bibitem{W}
Widodo,
\newblock Topological entropy of shift function on the sequences space induced by expanding piecewise linear transformations 
\newblock {\em Disc. and Cont. Dyn. Systems - Series A},
\newblock V. 8,No.1, 191--208, 2002


\end{thebibliography}
\end{document}